%% file: nfree_bruns.tex
\documentclass{amsart}
\usepackage{graphics}
\usepackage{amsfonts}
\usepackage{amssymb}
\usepackage[mathcal]{eucal}
\usepackage{mathrsfs}
\usepackage{latexsym}
\usepackage{amsthm}
\usepackage{amsmath,amscd}
\usepackage{enumerate}
\usepackage{epic}
\usepackage{hyperref}
\input{Macros}
\title[A simultaneous generalization of independence and disjointness]{A simultaneous generalization of independence and disjointness in boolean algebras}

\author[C.T. Bruns]{Corey~T. Bruns}
\address{Department of Mathematical and Computer Sciences, University of Wisconsin-Whitewater, 800 W. Main Street, Whitewater, WI 53190, USA}
\email{brunsc@uww.edu}
\urladdr{\url{http://facstaff.uww.edu/brunsc}}


\begin{document}
\begin{abstract}
We give a definition of some classes of boolean algebras generalizing free boolean algebras; they satisfy a universal property that certain functions extend to homomorphisms. We give a combinatorial property of generating sets of these algebras, which we call n-independent. The properties of these classes (n-free and $\omega$-free boolean algebras) are investigated. These include connections to hypergraph theory and cardinal invariants on these algebras. Related cardinal functions, $\mathfrak{i}_n$, the minimum size of a maximal n-independent subset and $\mathfrak{i}_\omega$, the minimum size of an $\omega$-independent subset, are introduced and investigated. The values of $\mathfrak i_n$ and $\mathfrak i_\omega$ on \PoF{} are shown to be independent of ZFC. 
\end{abstract}
\maketitle

\input{NFree}

\input{OtherBAs}

\bibliographystyle{spmpsci}
\bibliography{biblio}

\end{document}

%% file: Macros.tex
\long\def\symbolfootnote[#1]#2{\begingroup\def\thefootnote{\fnsymbol{footnote}}\footnote[#1]{#2}\endgroup}

\renewcommand{\P}{\mathscr{P}}

\newcommand{\G}{\mathscr{G}}

\newcommand{\T}{\ensuremath{\mathscr{T}}}
\renewcommand{\H}{\mathscr{H}}
\theoremstyle{plain}
\newtheorem{theorem}{Theorem}[section]
\newtheorem{lemma}[theorem]{Lemma}

\newtheorem{proposition}[theorem]{Proposition}
\newtheorem{corollary}[theorem]{Corollary}

\theoremstyle{definition}
\newtheorem{definition}[theorem]{Definition}
\newcommand{\leftexp}[2]{{\vphantom{#2}}^{#1}{#2}}
\newcommand{\sub}{\subseteq}
\newcommand{\sm}{\smallsetminus}

\newcommand{\rng}{\mathop{\mathrm{rng}}}

\newcommand{\fct}{\mathop{\longrightarrow}}

\newcommand{\Ult}{\mathop{\mathrm{Ult}}}
\newcommand{\Clop}{\mathop{\mathrm{Clop}}}
\newcommand{\cmpn}{\mathop{\mathrm{cmpn}}}
\newcommand{\res}{\upharpoonright}
\renewcommand{\phi}{\varphi}

\newcommand{\FR}[1]{\mathop{\mathrm{Fr}}\left(#1\right)}
\newcommand{\FinCo}[1]{\mathop{\mathrm{FinCo}}\left(#1\right)}
\newcommand{\At}[1]{\mathop{\mathrm{At}}\left(#1\right)}
\newcommand{\BC}[1]{\mathscr B_c\left(#1\right)}
\newcommand{\BA}[1]{\mathscr B_a\left(#1\right)}

\newcommand{\bas}{boolean algebras}
\newcommand{\ba}{boolean algebra}
\newcommand{\CSP}{countable separation property}
\newcommand{\PInd}{n\mathrm{Ind}}
\newcommand{\pa}[1]{\left(#1\right)}
\newcommand{\Bco}[1]{\mathscr B_{co}\pa{#1}}
\newcommand{\Baco}[1]{\mathscr B_{aco}\pa{#1}}

\newcommand{\J}{\mathfrak i}
\newcommand{\s}[1]{\left\{#1\right\}}
\newcommand{\pc}[1]{\ensuremath{\left(\perp #1\right)}}
\newcommand{\defeq}{\ensuremath{\stackrel{\text{\tiny def}}{=}}}
\renewcommand{\epsilon}{\varepsilon}
\renewcommand{\emptyset}{\varnothing}

\newcommand{\Wolog}{Without loss of generality}
\newcommand{\wolog}{without loss of generality}
\newcommand{\br}[1]{\left[#1\right]}
\renewcommand{\bar}{\overline}

\newcommand{\subs}[2]{\left[#1\right]^{#2}}

\newcommand{\PoF}{\ensuremath{\P\left(\omega\right)/{\mathrm {fin}}}}

%% file: NFree.tex
\section{Definitions}
A \ba{}  $A$ is free over its subset $X$ if it has the universal property that every function $f$ from $X$ to a \ba{}  $B$ extends to a unique homomorphism.  This is equivalent to requiring that $X$ be independent and generate $A$ (uniqueness).   A generalization, $\perp$-free, is introduced in Heindorf \cite{Heindorf94}, and some of its properties are dealt with.  I follow his notation for some of its properties, but that of Koppelberg \cite{Handbook} for the operations $+,\cdot,-,0,1$ on Boolean Algebras, with the addition that for an element $a$ of a \ba{}, we let $a^0=-a$ and $a^1=a$.  An elementary product of $X$ is an element of the form $\prod_{x\in R}x^{\epsilon_x}$ where $R$ is a finite subset of $X$ and $\epsilon\in \leftexp R2$.  We further generalize the notion of freeness to $n$-freeness for $1\leq n \leq \omega$.

It is nice to have a symbol for disjointness; we define $a\perp b $ if and only if $a\cdot b =0$.

\begin{definition}\label{def:npp} Let $n$ be a positive integer, $A$ and $B$ be nontrivial \bas{}, and $U\sub A$.  A function
$f:U\rightarrow B$ is $n$-preserving if and only if for every
$a_0,a_1,\ldots,a_{n-1}\in U$, $\prod_{i<n}a_i=0$ implies that
$\prod_{i<n}f\pa{a_i}=0$.
\end{definition}

An infinite version of this is also important.
\begin{definition}\label{def:opp} Let $A$ and $B$ be nontrivial \bas{}, and $U\sub A$.  A function
$f:U\rightarrow B$ is $\omega$-preserving if and only if for every
finite $H\sub U$, $\prod H=0$ implies that
$\prod f\left[H\right]=0$.
\end{definition}

Then we say that $A$ is $n$-free over $X$ if every $n$-preserving function from $X$ into arbitrary $B$ extends to a unique homomorphism.  The uniqueness just requires that $X$ be a generating set for $A$. 

The existence of such extensions is equivalent to an algebraic property of $X$, namely that $X^+$ is $n$-independent. This notion is defined below, and the equivalence is proved.  (For $n=1$, this is the usual notion of free and independent; for $n=2$, the notions are called $\perp$-free and $\perp$-independent by Heindorf \cite{Heindorf94};  Theorem 1.3 in the same paper shows that a $2$-free \ba{} has a $2$-independent generating set.  We differ from Heindorf in that he allows $0$ to be an element of a $\perp$-independent set.) Since any function that is $n$-preserving is also $m$-preserving for all $m\leq n\leq \omega$, so that an $m$-free \ba{} is also $n$-free over the same set; in particular, any $n$-free \ba{} is $\omega$-free. It's also worth noting that a function is $\omega$-preserving if and only if it's $n$-preserving for all finite $n$.

Freeness over $X$ implies that no elementary products over $X$ can be $0$.  $n$-independence weakens this by allowing products of $n$ or fewer elements of $X$ to be $0$.  This requires some other elementary products to be $0$ as well--if $x_1\cdot x_2\cdot\ldots\cdot x_m=0$, then any elementary product that includes $x_1,\ldots,x_m$ each with exponent $1$ must also be $0$.

\begin{definition} Let $A$ be a \ba{}.  For $n$ a positive integer, $X\sub A$ is $n$-independent if and only if $0\notin X$ and for all nonempty finite
subsets $F$ and $G$ of $X$, the following three conditions hold:
\begin{description}
\item[$\left(\perp 1\right)$] $\sum F \neq 1.$
\item[$\left(\perp 2\right)_n$] If $\prod F =0$, there is an $F'\sub F$ with
$\left|F'\right|\leq n$ such that $\prod F'=0$.
\item[$\left(\perp 3\right)$] If $0\neq \prod F \leq \sum G$, then $F\cap
G\neq\emptyset$.
\end{description}
\end{definition}

\begin{definition} Let $A$ be a \ba{}.  $X\sub A$ is $\omega$-independent if and only if $0\notin X$ and for all nonempty finite
subsets $F$ and $G$ of $X$, the following two conditions hold:
\begin{description}
\item[$\left(\perp 1\right)$] $\sum F \neq 1.$
\item[$\left(\perp 3\right)$] If $0\neq \prod F \leq \sum G$, then $F\cap
G\neq\emptyset$.
\end{description}
\end{definition}

We note that in both the above definitions, if $X$ is infinite, then
$\pc3\Rightarrow\pc1$; suppose \pc1 fails; take a finite $G$ with $\sum G=1$, then take some $x\notin G$
and let $F\defeq\s x$; then $0<\prod F\leq\sum G$ and $F\cap G\neq \emptyset$.

\pc3 has several equivalent forms which will be useful in the sequel.
\begin{proposition}\label{prop:a-2}The following are equivalent for a subset of $X$ of a \ba{}
A:
\begin{enumerate}
\item For all nonempty finite $F,G\sub X$, \pc3.
\item For all nonempty finite $F,G\sub X$ such that $F\cap G=\emptyset$ and
$\prod F\neq 0$, $\prod F \not\leq \sum G$.
\item For all nonempty finite $F,G\sub X$ such that $F\cap G=\emptyset$ and
$\prod F\neq 0$, $\prod F \cdot \prod -G\neq 0$, where $-G\defeq \s{-g: g\in
G}$.
\item Let $X$ be bijectively enumerated by $I$ such that $X=\s{x_i:i\in I}$. For all nonempty finite $R\sub I$ and all $\epsilon\in \leftexp R 2$ such
that $1\in\rng\epsilon$ and $\prod_{\stackrel{i\in R}{\epsilon_i = 1}}x_i\neq 0$,
$\prod_{i\in R}x_i^{\epsilon_i}\neq 0$. 
\end{enumerate}
\end{proposition}
In words, the final equivalent says that no elementary product of elements of $X$ is
$0$ unless the product of the non-complemented elements is $0$.  We note that
in the presence of $\pc2_n$, the words ``of $n$'' may be inserted after
``product.''
\begin{proof}
We begin by pointing out that \pc3 has two hypotheses, $0\neq \prod F$ and
$\prod F\leq\sum G$.  Thus the contrapositive of \pc3 is ``If $F\cap
G=\emptyset$, then $0=\prod F$ or $\prod F\not\leq\sum G$,'' which is
equivalent to (2).

(2) and (3) are equivalent by some elementary facts: $a\leq b \;\iff a\cdot-b=0$ and de Morgan's law that $-\sum G=\prod -G$.

(3) $\Rightarrow$ (4):

Assume (3) and the hypotheses of (4).  If $\rng \epsilon =\s1$, the conclusion is clear.  Otherwise, let $F\defeq\s{x_i:i\in R\mbox{ and } \epsilon_i=0}$.  Then (3) implies that $\prod_{i\in R}x_i^{\epsilon_i}\neq 0$, as we wanted.

(4) $\Rightarrow$ (3):

Assume (4) and the hypotheses of (3).  Let $R\defeq\s{i\in I : x_i\in F\cup G}$ and let $\epsilon_i=1$ if $x_i\in F$ and $\epsilon_i=0$ otherwise.  Then (4) implies that $\prod F\cdot\prod -G \neq 0$, as we wanted.
\end{proof}

\begin{proposition}\label{prop:a-1}The following are equivalent for a subset $X$ of a \ba{}
A:
\begin{enumerate}

\item $X$ is $\omega$-independent
\item Let $X$ be bijectively enumerated by $I$ such that $X=\s{x_i:i\in I}$. For all nonempty finite $R\sub I$ and all $\epsilon\in \leftexp R 2$ such
that $\prod_{\stackrel{i\in R}{\epsilon_i = 1}}x_i\neq 0$,
$\prod_{i\in R}x_i^{\epsilon_i}\neq 0$.
\end{enumerate}
\end{proposition}
\begin{proof}

The proof is similar to that of proposition \ref{prop:a-2}.  \pc1 is taken care of since products over an empty index set are taken to be $1$ by definition.
\end{proof}

In the same spirit, we have an equivalent definition of $n$-independent.
\begin{proposition}\label{prop:a}Let $n$ be a positive integer or $\omega$, $A$ a nontrivial \ba{} and $X\sub A^+$.  $X$ is $n$-independent if and only if for every $R\in\left[X\right]^{<\omega}$ and every $\epsilon\in\leftexp R2$, if $\prod_{x\in R}x^{\epsilon_x}=0$ then there is an $R'\sub R$ with $\left|R'\right|\leq n$ such that $\epsilon\left[R'\right]=\s1$ and $\prod R' =0$.
\end{proposition}
\begin{proof}
If $n=\omega$, this is part of proposition \ref{prop:a-1}.

Let $n$ be a positive integer, $A$ a \ba{}, and $X\sub A^+$.

We first show that $n$-independent sets have the indicated property.

Assume that $X$ is $n$-independent; take $R\in\left[X\right]^{<\omega}$ and $\epsilon\in\leftexp R2$ such that  $\prod_{x\in R}x^{\epsilon_x}=0$.  Let $F=\s{x\in R:\epsilon_x=1}$ and $G=\s{x\in R:\epsilon_x=0}$.  $F\neq\emptyset$; otherwise $\sum R=-\prod R = 1$, contradicting \pc1.  Since $\prod_{x\in R}x^{\epsilon_x}=\prod F \cdot \prod -G$, we have that $\prod F\leq\sum G$. 
If $G=\emptyset$, then $\sum G=0$ and so $\prod F =0$ as well.  If $G\neq \emptyset$, then $\prod F=0$ since $F\cap G=\emptyset$, using \pc3.  Then $R'$ is found by $\pc2_n$.
Now we show that sets with the indicated property are $n$-independent.

Assume that $X$ has the indicated condition and $F,G\in\left[X\right]^{<\omega}\sm\s\emptyset$.  We have three conditions to check.
\begin{description}
\item[\pc1]  Suppose that $\sum F=1$.  We let $F$ be the set $R$ in the condition, setting $\epsilon_x=0$ for all $x\in F$.  Then $\prod_{x\in F}x^{\epsilon_x}=\prod -F=-\sum F = 0$ and $\s{x\in F:\epsilon_x=1}=\emptyset$, thus there is no $R'$ as in the condition, since products over an empty index set are equal to 1.

\item[$\pc2_n$] Suppose that $\prod F=0$.  Again we let $F$ be the set $R$ in the condition, this time setting $\epsilon_x=1$ for all $x\in F$.  Then the condition gives us the necessary $F'$

\item[\pc3]  Suppose that $0\neq\prod F\leq\sum G$ and $F\cap G=\emptyset$.  Let $R=F\cup G$ and $\epsilon\in \leftexp R2$ be such that $\epsilon\left[F\right]=\s1$ and $\epsilon\left[G\right]=\s0$.  Then $\prod_{x\in R}x^{\epsilon_x}=0$ and the condition gives $\prod F=0$, which contradicts the original supposition.
\end{description}
\end{proof}


\begin{lemma}\label{pindtoind} If $H$ is an 
$\omega$-independent set that has no finite subset $F$ such that $\prod F=0$, $H$ is in fact
independent.  Furthermore, if $H$ is $n$-independent with no subset $F$ of size $n$ or less with $\prod F=0$, then $H$ is independent.\end{lemma}
\begin{proof}
We only need show that $\pc2_1$ holds, which it does vacuously.
\end{proof}

$2$-independence, and thus $n$-independence for $2\leq n\leq\omega$, is also a generalization of pairwise disjointness on
infinite sets.

\begin{theorem}\label{ipd}If $X\sub B^+$ is an infinite pairwise disjoint set, then $X$ is
$2$-independent.\end{theorem}
\begin{proof}

This is clear from proposition \ref{prop:a}.
\end{proof}

Some non-trivial examples of $2$-free \bas{} are the finite-cofinite algebras.  For infinite $\kappa$, let $A=\mathop{\mathrm FinCo}\left(\kappa\right)$.  $\mathrm{At}\pa A$  is a $2$-independent generating set for $A$.

Having an
$n$-independent generating set is equivalent to $n$-freeness. This is known in
Koppelberg \cite{Handbook} for $n=1$ and Heindorf \cite{Heindorf94} for $n=2$.  Our proof is more elementary than that of Heindorf \cite{Heindorf94} in that it avoids clone theory.

\begin{theorem}\label{thm:ofree_oindep}If $A$ is $\omega$-free over $X$, then
$X^+$ is $\omega$-independent.
\end{theorem}
\begin{proof}
Let $A$ and $X$ be as in the hypothesis; we show that $X^+$ is
$\omega$-independent.

\Wolog{}, we may assume that $0\notin X$ so that $X^+=X$.
\begin{description}
\item[\pc1]Let $f:X\rightarrow \s{0,1}$ be such that $f\br X =\s 0$.  Clearly $f$ is
$\omega$-preserving and thus extends to a homomorphism $\overline f$.  Take
$F\in\left[X\right]^{<\omega}$; then $\overline f \pa{\sum F}=\sum
f\left[F\right]=0$, so that $\sum F\neq 1$.
\item[\pc3]Take $F,G\in\left[X\right]^{<\omega}$ such that $F\cap G=\emptyset$ and
$\prod F\neq 0$.  Let $f:X\rightarrow\s{0,1}$ be such that
$f\left[F\right]=\s1$ and $f\left[X\sm F\right]=\s0$.  We claim that $f$ is
$\omega$-preserving.  If $H\sub X$ is finite such that
$\prod f\left[H\right]\neq0$, then it must be
that $H\sub F$, and hence $\prod H\neq 0$.  Thus $f$ extends to a homomorphism $\overline
f$.  Then $$\overline f\pa{\prod F\cdot \prod -G}=\prod
f\left[F\right]\cdot\prod\overline f\left[-G\right]=1,$$ and so $\prod
F\cdot\prod -G\neq 0$.
\end{description}

\end{proof}

\begin{theorem}\label{thm:nfree_nindep}  Let $n$ be a positive integer and $A$ a \ba{}.  If $A$ is $n$-free over $X$, then $X^+$ is $n$-independent.
\end{theorem}

\begin{proof} Again, \wolog{} $X=X^+$.

From theorem \ref{thm:ofree_oindep}, $X$ is $\omega$-independent, so we need
only show that $\pc2_n$ holds for $X$.  We do this by contradiction; assume
that $F\sub X$ is finite, of cardinality greater than $n$, $\prod F=0$, and
every subset $F'\sub F$ where $F'$ is of size $n$ is such that $\prod F'\neq
0$.

Define $f:X\rightarrow \s{0,1}$ by letting $f\left[F\right]=\s1$ and
$f\left[X\sm F\right]=\s0$.

Then $f$ is $n$-preserving. Let $G\sub X$ be of size $n$ and have $\prod G=0$.
Then $G\not\sub F$, so some $x\in G$ has $f\pa x =0$, so  $\prod f\left[G\right] = 0$.  Thus $f$ must extend to a homomorphism, but then $f\pa 0 =f\pa{\prod F}=\prod f\left[F\right]=\prod \s1=1$, which is a contradiction.
\end{proof}

\begin{theorem}\label{thm:oindep_ofree}Let $A$ be generated by its
$\omega$-independent subset $X$.  Then $A$ is $\omega$-free over $X$.
\end{theorem}
\begin{proof}
 Let $f$ be
an $\omega$-preserving function with domain $X$; we will show that
$f$ extends to a unique homomorphism.

Take a finite $H\sub X$ and $\epsilon\in\leftexp H2$ such that $\prod_{h\in
H}h^{\epsilon_h}=0$.  Then by \pc3 and \pc1, $\prod_{\epsilon_h=1}h=0$.  Then
since $f$ is $\omega$-preserving, $\prod_{\epsilon_h=1}f\pa h =0$ and
thus $\prod_{h\in H}f\pa h ^{\epsilon_h}=0$.  Thus by Sikorski's extension criterion, $f$ extends to a homomorphism. 

Uniqueness is clear as $X$ is a generating set.
\end{proof}

\begin{theorem}\label{thm:nindep_nfree}Let $n$ be a positive integer.  If $X$ generates $A$ and $X$ is $n$-independent, then $A$ is
$n$-free over $X$
\end{theorem}
\begin{proof}
 Let
$f$ be an $n$-preserving function with domain $X$.  

Take any distinct $x_0,x_1,\ldots,x_{k-1}$ and $\epsilon\in\leftexp k2$ such that
$\prod_{i<k}x_i^{\epsilon_i}=0$.  

Then by proposition \ref{prop:a}, there is an $F'\sub\s{x_i:\epsilon_i=1\;\mbox{and}\; i<k}$ such that $\left|F'\right|\leq
n$ and $\prod F'=0$.  Since $f$ is $n$-preserving, it must be that
$\prod f\left[F'\right]=0$, and thus $\prod_{i<k}f\pa{x_i}^{\epsilon_i}=0$.
Thus, by Sikorski's extension criterion, $f$ extends to a homomorphism. 

Uniqueness is clear as $X$ is a generating set.
\end{proof}

So we have shown that the universal algebraic property defining $n$-free \bas{}
is equivalent to having an $n$-independent generating set.

\begin{theorem}\label{semigroup} $\omega$-free \bas{} (and thus all $n$-free
\bas{}) are semigroup algebras.
\end{theorem}
A semigroup algebra is a \ba{}  that has a generating set that includes $\left\{0,1\right\}$, is closed under the product operation, and is disjunctive when $0$ is removed.
\begin{proof}
Let $A$ be $\omega$-free over $G$.  Then let $H'$ be the closure of $G\cup \left\{0,1\right\}$ under finite products, that is, the set of all finite products of elements of $G$, along with $0$ and $1$.  Clearly $H'$ generates $A$, includes $\left\{0,1\right\}$, and is closed under products, so all that remains is to show that $H=H'\sm\left\{0\right\}$ is disjunctive.
From proposition 2.1 of Monk \cite{CINV2}, $H$ is disjunctive if and only if for every $M\sub H$ there is a homomorphism $f$ from $\left<H\right>$ into $\P\left(M\right)$ such that $f\left(h\right)=M\downarrow h$ for all $h\in H$. 

To this end, given $M\sub H$, let $f:G\rightarrow \P\left(M\right)$ be defined by $g\mapsto M\downarrow g$.  We claim that $f$ is $\omega$-preserving.  Suppose $G'\in \left[G\right]^{<\omega}$ is such that $\prod G'=0$.  Then $\prod_{g\in G'}f\pa g = \bigcap_{g\in G'}\left(M\downarrow g \right)=\left\{a\in M :\forall g\in G'\left[a\leq g\right] \right\}=\emptyset$.  So $f$ extends to a unique homomorphism $\hat f$ from $A$ to $\P\left(M\right)$.  If $h\in H\sm \left\{1\right\}$, then $h=g_1 \cdot g_2 \cdot \ldots \cdot g_n$ where each $g_i\in G$.  So $$\hat f\left(h\right)=\hat f\left(g_1 \cdot g_2 \cdot \ldots \cdot g_n\right)=f\left(g_1\right)\cap f\left(g_2\right)\cap \ldots \cap f\left(g_n\right)=$$ $$\left(M\downarrow g_1 \right)\cap \left(M \downarrow g_2\right)\cap \ldots \cap \left(M\downarrow g_n\right) =M\downarrow \left(g_1 \cdot g_2 \cdot \ldots \cdot g_n\right)=M\downarrow h.$$  Likewise, $\hat f\left(1\right)=M = M\downarrow 1$.  Thus $H$ is disjunctive and $A$ is a semigroup algebra over $H'$.
\end{proof}

\section{Hypergraphs and their Anticlique Algebras}
There is a correspondence with hypergraphs for $\omega$-free \bas{}.
We recall that a hypergraph is a pair $\G=\left<V,E\right>$ where $V$ is called the vertex set, and $E\sub\P\pa V\sm\s\emptyset$ is called the hyperedge set; an element of $E$ is called a hyperedge. We will insist on loopless hypergraphs, that is, $E\sub\P\pa V\sm\subs{V}{\leq 1}$.
 A hypergraph is $n$-uniform if $E\sub\left[V\right]^n$.  For a given hypergraph, we call a set $A\sub V$ an anticlique if it includes no hyperedges; that is, for all $e\in E$, $e\sm A\neq \emptyset$, and call the set of anticliques $A\pa \G$. 
Given a hypergraph $\G$, we define an $\omega$-free \ba{} as a subalgebra of $\P\pa{A\pa \G}$.  For $v\in V$, let $v_+\defeq\s{A\in A\pa\G : v\in A}$, which is an element of $P\pa{A\pa \G}$, and for a set $H$ of vertices, $H_+\defeq\s{v_+:v\in H}$.  We then define the anticlique algebra of $\G$ as $\BA\G\defeq\left<V_+\right>$.

We do not consider cliques in general hypergraphs; it's not clear which way to define them.  For an $n$-uniform hypergraph, a clique may be non-controversially defined as a set $C$ where $\subs C n \sub E$, but for a  hypergraph with hyperedges of different cardinalities, it is not clear how many hyperedges must be included in a clique.    This difficulty stems from a lack of a reasonable way to define ``complement hypergraph.''  A few possibilities for the hyperedge set of $\bar\G$ are $\P\pa G \sm E$,  $\subs G {<\omega} \sm E$,  and $\subs G {\leq\pa{\sup_{e\in E}\pa{\left|e\right|}}}\sm E$.  For an $n$-uniform hypergraph $\left<G,E\right>$, the complementary hypergraph is $\left<G,\subs G n \sm E\right>$, and then a clique in $\G$ is an anticlique in $\bar\G$.  Each possible definition for complement hypergraph results in a different definition for clique, all of which are more complicated than our definition of anticlique. Since anticliques suffice for our study, we do not choose a side on what a clique ought to be.

\begin{theorem}For any hypergraph $\G=\left<V,E\right>$, $\BA\G$ is $\omega$-free over $V_+$.\end{theorem}
\begin{proof} We need only show that $V_+$ is $\omega$-independent; we will use proposition \ref{prop:a}.

Suppose that $R\in\subs V{<\omega}, \epsilon\in \leftexp R 2$, and $\bigcap_{v\in R} v_+^{\epsilon_v}=\emptyset$.  Let $S=\s{v\in R:\epsilon_v=1}$.  If $\bigcap_{v\in S} v_+\neq \emptyset$, let $T$ be a member of $\bigcap_{v\in S} v_+$.  So then $T$ is an anticlique, and $S\sub T$.  We note that clearly every subset of an anticlique is again an anticlique, so $S$ is also an anticlique, and $S\in\bigcap_{v\in R}v_+^{\epsilon_v}$.  

\end{proof}

If the hypergraph is somewhat special, we have more:
\begin{theorem}\label{thm:hyperedgesize}For any hypergraph $\G=\left<V,E\right>$ where $E\sub\left[V\right]^{\leq n}$, $\BA\G$ is $n$-free.\end{theorem}
\begin{proof}We show that $V_+$ is $n$-independent.

From the previous theorem, we need only show that $\pc2_n$ holds for $V_+$.  Let $F$ be a finite subset of $V$ such that $\prod F_+=0$.  Using the observation that $\prod F_+$ is the set of anticliques that include $F$, $F$ is not an anticlique.  Thus some hyperedge $e$ is a subset of $F$.  Then $\prod e_+=0$ as no anticlique can include that hyperedge.  Since all hyperedges have at most $n$ vertices, $\left|e\right|\leq n$, which is what we wanted.  \end{proof}

We also reverse this construction. Given a \ba{} $A$ with an $\omega$-independent generating set $H$, we construct a hypergraph $\G$ such that $A\cong\BA\G$; we call it the $\perp$-hypergraph of $A,H$.  The vertex set is $H$, and the hyperedge set is defined as follows; a subset $e$ of $H$ is a hyperedge if and only if the following three conditions are all true:
\begin{enumerate}
\item $e$ is finite.
\item $\prod e=0$.
\item If $f\subsetneq e$, then $\prod f\neq 0$.
\end{enumerate}

We have only finite hyperedges in this graph, and no hyperedge is contained in another.  Note that if $H$ is $n$-independent, the hyperedge set is included in $\left[H\right]^{\leq n}$.

\begin{theorem}\label{thm:hypergraphofalg}
Let $n$ be a positive integer or $\omega$, $X\sub A$ be $n$-independent and generate $A$, and $\G=\left<X,E\right>$ be the
$\perp$-hypergraph of $A$.  Then $A\cong\BA\G$.\end{theorem}
\begin{proof}Let $f:X\rightarrow X_+$ be defined so that $v\mapsto v_+$ for $v\in X$.  We claim that $f$ is an
$n$-preserving function.  If $G\sub X$ is of size $\leq n$ such that
$\prod G =0$, then it has a subset $G'$ minimal for the property of having $0$
product; thus $G'\in E$, so that $\prod G'_+=0$, and so $\prod
f\left[G\right]=0$.

$f$ is bijective, and its inverse is also $n$-preserving; the image of
$f$ is a generating set, so that $f$ extends to an isomorphism.
\end{proof}

\begin{definition}Let $\G_i\left<V_i,E_i\right>$ be hypergraphs for $i\in\s{0,1}$.  A hypergraph homomorphism is a function $f:V_0\rightarrow V_1$ such that if $e\in E_0$, then $f\left[e\right]\in E_1$.
\end{definition}

Notice that a graph homomorphism is a hypergraph homomorphism when the graphs are considered as 2-uniform hypergraphs.



In the rest of this section we consider ordinary graphs, that is, hypergraphs for which $E\sub V 2$. 
In this case, ``clique'' is not ambiguous, so we can define the clique algebra of a graph. We let $C\pa\G$ be the set of cliques in $\G$, and $v_+$ be the set of cliques including vertex $v$. (This conflicts with an earlier use of $v_+$, but context will make it clear which is meant.) Then $\BC\G$ is the subalgebra of $\P\left(C\left(\G\right)\right)$ generated by  $\left\{v_+ : v\in
G\right\}$.

%
%
%
%
%
%


We give some examples of $2$-free \bas{} with unusual properties.

For a $2$-free algebra of the form $\BC T$ for a tree (in the graph-theoretical
sense--a connected acyclic graph) or a forest $T$ of size
$\kappa$, there
are further conclusions that can be drawn.  As a forest has no triangles, all the cliques in $T$ are of size at most 2.

So any subset
of $T_+$  of size 3 or more has a disjoint pair.  

If $T$ is a $\kappa$-tree (in the order theoretic sense, 
that is, of height $\kappa$ and each level of size $<\kappa$),and we take the edge set to consist of pairs $\s{u,v}$ where $v$ is an immediate successor of $u$, then $T_+$ has
a pairwise disjoint subset of size $\kappa$--take an element of every other
level--so that $\FinCo\kappa\leq\BC T$, and $\FR\kappa\leq \BA T$. 

 It seems
difficult to avoid one of $\FinCo\kappa$ and $\FR\kappa$ as a subalgebra, as it
is necessary to find a graph of size $\kappa$ with no clique or anticlique of
size $\kappa$.  A witness to $\kappa\not\longrightarrow \pa\kappa^2_2$ is the edge set of such a graph, but we do not know about the variety of such witnesses.  If $\kappa$ is weakly compact, then there are no such witnesses and so for any graph of size $\kappa$, $\FinCo\kappa$ or $\FR\kappa$ is a subalgebra of $\BC\G$ 

As a graph can be characterized as a symmetric non-reflexive relation, for any non-reflexive relation $R$,
we may form  algebras $\BA{R\cup R^{-1}}$ and $\BC{R\cup R^{-1}}$.  When $R$
is an ordering of some sort, $R\cup R^{-1}$ is usually called the (edge set of
the) comparability
graph of $R$.  Thus for a (non-reflexive
) ordering $\left<P,<\right>$, it has comparability
graph $\G_P=\left<P,<\cup<^{-1}\right>$ and  we define its
comparability algebra $\Bco P \defeq\BC{\G_P}$ and its
incomparability algebra $\Baco P \defeq\BA{\G_P}$.  Since points in the partial order are vertices of the comparability graph, we may use the $p_+$ notation without fear of confusion. 
When $P$ is a partial order in the strict sense, $C\sub P$ is a clique in
$\G_P$ if and only if $C$ is a chain in $\leq$ if and only if $C_+$ is an independent subset of
$\Bco P$, and $A\sub P$ is an anticlique in
$\G_P$ if and only if $A$ is an antichain in $\leq$ if and only if $A_+$ is a pairwise disjoint
set in $\Bco P$.
So if $\left<T,\leq\right>$ is a $\kappa$-Suslin tree, in both $\Bco T$ and $\Baco T$, $T_+$ is a
$2$-independent set of size
$\kappa$, but has no independent subset of size $\kappa$, nor a pairwise
disjoint subset of size $\kappa$ since $T$ has neither chains nor antichains of
size $\kappa$.

\begin{proposition}
If $f:P\rightarrow Q$ is a strictly order-preserving function, that is, a morphism in the category of partial orders, then there is a homomorphism $f^*:\Baco P \rightarrow \Baco Q$ such that $f^*\pa{p_+}=f\pa p _+$.
\end{proposition}
\begin{proof}
By the universal property of $2$-free \bas{}, we need only show that $g$ is $2$-preserving where $g\pa{p_+}=f\pa p_+$; then $g$ extends to the $f^*$ of the conclusion.

Fix distinct $p,p'\in P$; if $p_+ \perp p'_+$ in $\Baco P$, then $p$ and $p'$ are comparable in $P$, without loss of generality, $p<p'$.  Then $f\pa p < f\pa{ p'}$, so that $f\pa p_+ \perp f \pa {p'}_+$.
\end{proof}
Similarly, an incomparability-preserving map from $P$ to $Q$ gives rise to a homomorphism of $\Bco P$ and $\Bco Q$.

\section{Hypergraph Spaces}

The dual spaces to $\omega$-free \bas{} are also interesting.  Like with
graphs, a hypergraph space may be defined in terms of a hypergraph--the
definition generalizes that of a graph space.

\begin{definition}
Let $\G=\left<G,E\right>$ be a hypergraph and  $A\pa\G$ its set of anticliques.  For each $v\in
G$, we define two sets: $$v_+\defeq\s{A\in A\pa\G : v\in A}$$  $$v_-\defeq\s{A\in A\pa\G
: v\notin A}.$$

Then the hypergraph space of $\G$ is the topology on $A\pa\G$ with $\bigcup_{v\in G} \s{v_+,v_-}$ as a closed subbase.

Any topological space $\T$ for which there is a hypergraph $\G$ such that $\T$ is homeomorphic to the hypergraph space of $\G$ is called a hypergraph space.

\end{definition}

\begin{theorem} The Stone dual of an $\omega$-free \ba{} is a hypergraph space.
\end{theorem}
\begin{proof}  Let $A$ be an $\omega$-free \ba{}.  Thus by theorem \ref{thm:hypergraphofalg}, there is a hypergraph $\G$ such that $A\cong \BA\G$.  Let $\T$ be the hypergraph space of $\G$.  We claim that $\Clop\pa\T\cong A$.

In fact, $\Clop\pa\T = \BA\G$.  On both sides here, elements are sets of anticliques of $\G$.  As $\T$ is defined by a clopen subbase, elements of $\Clop\pa\T$ are finite unions of finite intersections of elements of that subbase $\bigcup_{v\in G} \s{v_+,v_-}$.  Elements of the right hand side are sums of elementary products of elements of $\bigcup_{v\in G} \s v_+$, that is, sums of finite products of elements of $\bigcup_{v\in G} \s{v_+,v_-}$.  As the operations are the usual set-theoretic ones on both sides, they are in fact the same algebra.  

The topological result follows by duality.
\end{proof}

We repeat a few definitions from Bell and van Mill \cite{BellvanMill80} needed for some topological applications.

\begin{definition}Let $n\in \omega$ for all these definitions.

A set $S$ is $n$-linked if every $X\in\left[S\right]^n$ has
non-empty intersection.

A set $P$ is $n$-ary if every $n$-linked subset of $P$ has non-empty
intersection.

A compact topological space $\T$ has compactness number at most $n$, written $\cmpn\pa\T \leq n$,
if and only if it has an $n$-ary closed subbase.
$\T$ has compactness number $n$, written $\cmpn\pa\T = n$, if and only if $n$ is the least integer
for which $\cmpn \pa \T \leq n$. $\cmpn\pa\T = \omega$ if there is no such $n$.
\end{definition}

The following generalizes and algebraizes proposition 3.1 of Bell \cite{Bell82}.

\begin{proposition}\label{prop:bell}If a \ba{} $A$ is $n$-free for some $2\leq n \leq \omega$, then $\cmpn\pa{\Ult A}\leq n$.\end{proposition}
\begin{proof} This is vacuous if $n=\omega$.  If $n<\omega$, then $\Ult\pa A$ is a hypergraph space for a hypergraph $\G$ with all hyperedges of size $\leq n$.

We take the clopen subbase $S=\bigcup_{v\in G} \s{v_+,v_-}$ of the hypergraph space of $\G$ and show that it is $n$-ary.  Let $\mathcal F\sub S$ be $n$-linked.  We may write $\mathcal F=\s{v_+ : v\in A}\cup \s{v_- : v\in B}$ for some $A,B\sub G$.  Since $v_+\cap v_-=\emptyset$ and $n\geq 2$, $A\cap B=\emptyset$.  Let $A'$ be a finite subset of $A$.  Since any product of $n$ or fewer elements of $\mathcal F$ is non-zero, $A'$ must be an anticlique in $\G$; if not, then $\prod A'_+=0$, so then $A'_+$ would have a subset of size $n$ with empty intersection, contradicting that $\mathcal F$ is $n$-linked.
Thus $ A'\in \bigcap \mathcal F$, that is, $\mathcal F$ has non-empty intersection and thus $S$ is $n$-ary.
\end{proof}

Bell's \cite{MR803933} corollary 5.2 shows that certain topologies on
$\left[\omega_1\right]^{\leq m}$ have compactness number $n$ for certain 
$n,m\leq \omega$.  These topologies are the hypergraph spaces of
$\left<\omega_1,\left[\omega_1\right]^{2n-3}\right>$ and
$\left<\omega_1,\left[\omega_1\right]^{2n-2}\right>$.
\begin{theorem}For infinitely many $n\in\omega$, there is a \ba{} which is $n$-free and is not $\pa{n-1}$-free.
\end{theorem}
\begin{proof}Let $k$ be the least integer for which $\BA{\left<\omega_1,\left[\omega_1\right]^{2n-3}\right>}$ is $k$-free and $\ell$ be the least integer for which $\BA{\left<\omega_1,\left[\omega_1\right]^{2n-2}\right>}$ is $\ell$-free. We have that
$n\leq k \leq 2n-3$ and $n\leq \ell\leq 2n-2$.  The lower bounds are a
consequence of the compactness numbers of those spaces (Bell's \cite{MR803933} result and proposition \ref{prop:bell}), while the upper bounds
are a consequence of theorem \ref{thm:hyperedgesize}.

Thus we have, for arbitrary $n\in\omega$, an $\omega$-free \ba{} of finite freeness at least
$n$.
\end{proof}

\section{Constructions}\label{sec:ac}

\newcommand{\GS}{$\mathscr{GS}$ }
In this section, we consider the categories of $n$-independently generated \bas{} and of hypergraph spaces and their behavior under some constructions.

If a \ba{} $A$ is $\omega$-free, it is isomorphic to $\BA{\G}$ for some hypergraph $\G$; we'll call this the $\perp$-hypergraph of $A$.  If a \ba{} is $2$-free, this $\perp$-hypergraph is a graph, so we can just call it the $\perp$-graph.  Such a \ba{} is also isomorphic to $\BC{\G}$ for for a graph $\G$, which is called the intersection graph of $A$.
 
We will show in section \ref{sec:cf} that complete \bas{} are not $\omega$-free.
As $\P\left(\kappa\right)$ is isomorphic to $\leftexp \kappa 2$, the class of $\omega$-free \bas{} is not closed under infinite products.

\begin{theorem}\label{finprod}Let $2\leq n\leq \omega$. If $H\sub A$ and $K\sub B$ are $n$-independent, then  
$L\defeq\left(H\times \left\{0\right\}\right) \cup \left(\left\{0\right\}\times K\right)$ is $n$-independent in $A\times B$.
\end{theorem}
\begin{proof}

We will apply proposition \ref{prop:a}.  Suppose that $F\in\subs H {<\omega}, G\in\subs K {<\omega}, \epsilon\in \leftexp F2, \delta\in\leftexp G2$, and $\prod_{x\in F}\pa{x,0}^{\epsilon_x}\cdot\prod_{y\in G}\pa{0,y}^{\delta_y}=0$.  If there are $x\in F$ and $y\in G$ such that $\epsilon_x=\delta_y=1$, then $\pa{x,0}\cdot\pa{0,y}=0$ as desired.  Otherwise, without loss of generality, we may assume that $\epsilon\left[F\right]\sub\s0$.  Then $\prod_{x\in F}\pa{x,0}^{\epsilon_x}=\pa{\prod_{x\in F}-x,1}$, so that $\prod_{y\in Y}y^{\delta_y}=0$; then the $n$-independence of $K$ gives the result.
\end{proof}

It is important to note that $L$ does not generate $\left<H\right>\times \left<K\right>$; in fact (theorem \ref{thm:fcf_not2free}), the product of $n$-free \bas{} is not in general $n$-free. However, it is
the case that $\left<H\right>\times \left<K\right>$ is a simple extension of the subalgebra generated by $L$; $\left<L\right>\pa{\pa{1,0}}=\left<H\right>\times \left<K\right>$.

This result generalizes to infinite products quite easily, though the notation
is considerably more cumbersome.

\begin{theorem}\label{bigbindset} For $2\leq n \leq \omega$, if $\left<A_i : i\in I\right>$ is a system of \bas{} and for every
$i\in I$, $H_i\sub
A_i$ is $n$-independent in $A_i$, then the set $H\defeq\bigcup_{i\in I}
p_i\left[H_i\right]$, where
$$p_i\left(h\right)\left(j\right)\defeq\left\{\begin{array}{ll}h & i=j \\ 0 & i\neq
j\end{array}\right.$$ is $n$-independent in $A\defeq\prod_{i\in I}A_i$ and
$\prod^{\mathrm w}_{i\in I}A_i$.
\end{theorem}
\begin{proof}

This is essentially the same as  theorem \ref{finprod} with
more cumbersome notation.  

$p_i\left(h\right)$ is the function in
$A$ that is $0$ in all but the $i$th coordinate and is $h$ in the $i$th
coordinate, so that the projections $\pi_i\left[p_i\left[H_i\right]\right]=H_i$ and for $i\neq
j$,  $\pi_j\left[p_i\left[H_i\right]\right]=\left\{0\right\}$.

We apply proposition \ref{prop:a}.  Suppose that $R\in \subs H {<\omega}, \epsilon\in\leftexp R2$, and $\prod_{x\in R}x^{\epsilon_x}=0$.  Let $J\defeq\s{i\in I:R\cap p_i\left[H_i\right]\neq \emptyset}$.  If $J$ is a singleton, say $J=\s i$, then the $n$-independence of $H_i$ clearly makes $H$ $n$-independent. So we now concern ourselves with the case that  $\left|J\right|>1$, that is, we have distinct $i,j\in J$. If there are $x\in p_i\left[H_i\right]$ and $y\in p_j\br{H_j}$ with $\epsilon_x=\epsilon_y=1$, then $x\cdot y=0$ and we have our conclusion.  So we may assume that there is at most one $i\in J$ for which there is an $x\in p_i\br{H_i}$ such that $\epsilon_x=1$.  Then for any particular $i\in I$, $\prod\s{x^{\epsilon_x}:x\in R, x\notin p_i\br{H_i}}$ has $i$-th coordinate $1$, and so the facts that 
$$0=\prod_{x\in R}x^{\epsilon_x} = \prod\s{x^{\epsilon_x}:x\in R\cap p_i\br{H_i}},$$
 and that all the $H_i$ are  $n$-independent make $H$ $n$-independent.

\end{proof}

When $n=2$, we can also consider the $\perp$-graph and intersection graph of $H$ in the above theorem.  
The intersection graph is easily described: two elements of $H$ have non-zero
product if and only if they have non-zero product in one of the factors, so that the
intersection graph is the disjoint union of the intersection graphs of the
$H_i$.  The $\perp$-graph is more complex.  The $\perp$-graph of each $H_i$ is
an induced subgraph, but these subgraphs are connected to each other--each
vertex in $H_i$ is connected to every vertex in $H_j$ for $i\neq j$. This
construction is the ``join".

In other words, for any collection $\left<\G_i\right>$ of graphs,  $\BC{\bigcup_{i\in I}\G_i}\leq \prod_{i\in I}\BC{\G_i}$ and 
$\BA{\biguplus_{i\in I}\G_i}\leq \prod_{i\in I}\BA{\G_i}$.

The use of the word ``free" in $n$-free is warranted by the following:

\begin{theorem}\label{freeprod}  Suppose that $A\defeq\bigoplus\limits_{\stackrel{C}{i\in I}} A_i$ is an amalgamated free product of subalgebras $A_i$ for $i\in I$, where $C\leq A_i$ for each $i\in I$, $A_i\cap A_j =C$ for $i\neq j$, $A_i$ is $n$-free over $H_i$, and $C\leq\left<H_i\cap H_j\right>$.  Then $A$ is $n$-free over $\bigcup_{i\in I}H_i$.

\end{theorem}
\begin{proof}For convenience, assume that each $A_i \leq A$, $C\leq A_i$ and that, for
$i\neq j$, $A_i\cap A_j = C$, and
that $H_i$ is a set over which $A_i$ is $n$-free.  
We show that $A$ is $n$-free over $H\defeq\bigcup_{i\in I} H_i$.  

Let $B$ be a \ba{}, and  $f:H\rightarrow B$ be $n$-preserving.  Then
for each $i\in I$, $f_i:=f\res H_i$ is also $n$-preserving.  So each $f_i$ extends to a unique homomorphism $\phi_i : A_i\rightarrow B$.  That $\phi_i\res C = \phi_j\res C$ is clear as $C\sub\left<H_i\cap H_j\right>$.

 Then the universal property of amalgamated free products gives a unique homomorphism $\phi : A\rightarrow B$ that extends every $\phi_i$.  Note that
$$\phi\res H = \phi \res\bigcup_{i\in I} H_i=\bigcup_{i\in I}\left(\phi\res H_i\right)=\bigcup_{i\in I} \left(\phi_i\res H_i\right) = \bigcup_{i\in I} f_i = f.$$
So we have a unique extension of $f$ to a homomorphism, which is what we
wanted.\end{proof}

This of course includes free products.

An example where $C\neq\s{0,1}$ is as follows:
Let $\G$ be the complete graph on the ordinal $\omega_1 +\omega$ and $\H$ the complete graph on the ordinal interval $\left(\omega_1,\omega_1\cdot 2\right)$.  Then $\BA\G\cong\BA\H\cong\FR{\omega_1}$.  Note that $G\cap H=\left(\omega_1,\omega_1+\omega\right)$ so that $G_+\cap H_+=\left(\omega_1,\omega_1+\omega\right)_+$; we let $C=\left<\left(\omega_1,\omega_1+\omega\right)_+\right>\cong\FR\omega$.  It is clear that $C$ is as required in theorem \ref{freeprod}. Then we have that $\BA\G\mathop{\oplus}\limits_C \BA\H$ is $2$-free over $G_+\cup H_+$.

If $C$ is $2$-free over $\bigcup_{i\in I}H_i$, the $\perp$-graph of $\bigcup_{i\in I}H_i$ is easily described in terms of those of
$H_i$. It
is the ``amalgamated free product'' or  ``amalgamated disjoint union'' in the category of graphs---the same universal property holds.  More concretely, given a set of graphs $\G_i=\left<G_i,E_i\right>$, each of which has $\mathscr F=\left<F,E\right>$ as a  subgraph, the amalgamated disjoint union of the $\G_i$ over $\mathscr F$ is a graph on the union of the vertex sets where two vertices are adjacent if and only if they are adjacent in some $\G_i$.  That is, elements of $G_i\sm F$ and $G_j\sm F$ are not adjacent for $i\neq j$.  

In case $C=2$ and we have a free product, the $A_i$ form a family of independent subalgebras, so two elements
of $H$ (constructed in the proof above) have product zero if and only if they are in the same $H_i$ and have zero product
in $A_i$.  So the $\perp$-graph of $H$ is the disjoint union of the $\perp$
graphs of the $H_i$.
The intersection graph of $H$ is similarly constructed from those of the $H_i$:
the independence of the $A_i$ means that the intersection graph of $H$ is the
join of the intersection graphs of the $H_i$.

That is, $\bigoplus_{i\in I}\BA{\G_i} = \BA{\bigcup_{i\in I} \G_i}$ and 
$\bigoplus_{i\in I}\BC{\G_i} = \BC{\biguplus_{i\in I}\G_i}$.

Products of $n$-free \bas{} behave in a somewhat more complicated manner.  As discussed previously, infinite products
of $\omega$-free \bas{} are not necessarily $\omega$-free.  

\begin{theorem}\label{thm:fcf_not2free}$\FinCo{\omega_1}\times\FR{\omega_1}$ is not $2$-free.\end{theorem} 
\begin{proof} We use subscript function notation for the coordinates of tuples; i.e. $\pa{a,b}_0=a$ and $\pa{a,b}_1=b$.  We also extend this to sets of tuples; $\s{\pa{a,b},\pa{c,d}}_0=\s{a,c}$.

  We proceed by contradiction; suppose that $A\defeq\FinCo{\omega_1}\times\FR{\omega_1}$ is $2$-free over $X$, where $0\notin X$, that is, $X$ is $2$-independent. 

Consider $a_\alpha\defeq\pa{\s\alpha,0}$ for $\alpha<\omega_1$.  $a_\alpha$ is
an atom in $A$, so it must be an elementary product of $X$, that is,
$a_\alpha=\prod_{x\in H_\alpha}x^{\epsilon\pa{\alpha,x}}$, with
$H_\alpha\in\left[X\right]^{<\omega}$.  So let
$M\in\left[\omega_1\right]^{\omega_1}$ be such that $\s{H_\alpha:\alpha\in M}$
is a $\Delta$-system with root $F$.  Let $G_\alpha\defeq H_\alpha\sm F$.
Since $$M=\bigcup_{\delta\in \leftexp F2}\s{\alpha\in M:
\forall x\in F\left[\epsilon\pa{\alpha,x}=\delta_x\right]},$$ there is an uncountable $N\sub M$
such that $\epsilon\pa{\alpha,x}=\epsilon\pa{\beta,x}$ for all $\alpha,\beta\in
N$ and all $x\in F$, so that we may write, for $\alpha\in N$, $a_\alpha=\prod_{x\in
F}x^{\delta_x}\cdot \prod_{x\in G_\alpha}x^{\epsilon\pa{\alpha,x}}$. 
For each $\alpha\in N$, let $G'_\alpha\defeq\s{x\in G_\alpha :
\epsilon\pa{\alpha,x}=1}$.  If $\alpha,\beta\in N$ with $\alpha\neq\beta$, then
there are $x\in G'_\alpha$ and $y\in G'_\beta$ such that $x\cdot y=0$, by the
$2$-independence of $X$ and the fact that $$0=a_\alpha\cdot a_\beta=\prod_{x\in F}x^{\delta_x} \cdot\prod_{x\in G_\alpha}x^{\epsilon\pa{\alpha,x}}\cdot\prod_{x\in G_\beta}x^{\epsilon\pa{\beta,x}},$$  thus $\prod G'_\alpha \cdot \prod G'_\beta=0$.   Since $\FR{\omega_1}$ has cellularity $\omega$, the set $\s{\alpha\in N:\pa{\prod G'_\alpha}_1\neq 0}$ is countable, hence $P\defeq N\sm \s{\alpha\in N:\pa{\prod G'_\alpha}_1\neq 0}$ is uncountable and for $\alpha\in P$, $\pa{\prod G'_\alpha}_1=0$.  Since $\pa{\prod G'_\alpha}_0\cdot\pa{\prod G'_\beta}_0=0$ for distinct $\alpha,\beta\in P$, each $\pa{\prod G'_\alpha}_0$ is finite when $\alpha\in P$.

$X$ must generate $\pa{1,0}$; let $b_j$ for $j<n$ be disjoint elementary
products of $X$ such that $\sum_{j<n}b_j=\pa{1,0}$.  Thus there must be exactly
one $i<n$ such that $b_{i0}$ is cofinite; without loss of generality, $i=0$ so that $b_{00}$ is cofinite
and $b_{01}=0$.  Write $b_0$ as an elementary product, that is
$b_0=\prod_{j<m}c_j^{\xi_j}$ with each $c_j\in X$.  Then choose an $\alpha\in
P$ such that $\prod G'_\alpha\leq b_0$ and $G'_\alpha\cap\s{c_j:j<m}=\emptyset$.
Then $\prod G'_\alpha\cdot\sum_{j<n}c_j^{1-\xi_j}=0$, so $\rng\xi=\s0$; that
is, $b_0=\prod_{j<m} -c_j$.

Note that $X_1$ generates $\FR{\omega_1}$, so it must be uncountable, thus  
\\$\left(X\sm\s{c_j:j<m}\sm F\right)_1$ is also uncountable; let $Y\sub X$ be such that $Y_1$ is
an uncountable independent subset of $\pa{X\sm\s{c_j:j<m}\sm F}_1;$ such a $Y$
exists by theorem 9.16 of Koppelberg \cite{Handbook}.  Note that no finite
product of elements of $Y$ is $0$. Let
$\theta:Y\rightarrow\s{0,1}$ be such that $d_y\defeq\pa{y^{\theta_y}}_0$ is finite for
each $y\in Y$.

Consider $\s{d_y:y\in Y}$; Each $d_y$ is finite and $Y$ is an uncountable set, and thus there is 
an uncountable $Z\sub Y$ where $\s{d_y:y\in Z}$ is a $\Delta$-system with root $r$.  Let $y,z,t\in Z$ be 
distinct.  Then let $e_y\defeq d_y\sm r$, $e_z\defeq d_z\sm r$, and $e_t\defeq d_t\sm r$.  Then $$d_y\cdot d_z \cdot -d_t=\pa{e_y \cup r}\cap \pa{e_z\cup r}\cap\pa{\omega_1\sm\pa{e_t\cup r}}=r\cap \pa{\omega_1\sm\pa{e_t\cup r}}=\emptyset,$$  Then $\prod_{j<n}-c_j\cdot y^{\theta_y}\cdot z^{\theta_z}\cdot t^{1-\theta_t}=0$ and again, the only elements with exponent $1$ are elements of $Y$ and thus there is no disjoint pair, contradicting proposition \ref{prop:a}.

So we have a contradiction and thus there is no $2$-independent generating set
for $A$.
\end{proof}
 This is also an example of a simple extension of a
$2$-free \ba{} that is not $2$-free; the full product is a simple extension by
$\left(0,1\right)$ of the subalgebra generated by the set in theorem
\ref{finprod}.

The dual of this theorem is that we have two graph spaces whose disjoint union
is not a graph space; in fact we can say a bit more since the disjoint union of two
supercompact spaces is supercompact.  We show a slightly more general result
here:

\begin{proposition}\label{prop:union-scomp}If $X$ and $Y$ are $n$-compact
spaces, then $X\dot\cup Y$ is $n$-compact.\end{proposition}
\begin{proof}Suppose that $S$ and $T$ are $n$-ary subbases for the closed sets of $X$ and $Y$
respectively; that is, for any $S'\sub S$ with $\bigcap S'=\emptyset$, there
are $n$ members $a_1,a_2,\ldots,a_n$ of $S'$ such that $a_1\cap
a_2\cap\ldots\cap a_n=\emptyset$, and similarly for $T$.  Then $W\defeq S\cup
T\cup\s{X,Y}$ is an $n$-ary subbase for the closed sets of $X\dot\cup Y$.\end{proof}

So, letting $n=2$, the dual space of $\FinCo{\omega_1}\times\FR{\omega_1}$ is supercompact, but
is not a graph space. 

\begin{theorem}$\FinCo{\omega_1}\times\FR{\omega_1}$ is $3$-free.\end{theorem}
\begin{proof} Let $\s{x_\alpha:\alpha<\omega_1}$ be an independent generating
set for $\FR{\omega_1}$.  Then the set
$X\defeq\s{\pa{\s\alpha,x_\alpha}:\alpha<\omega_1}\cup\s{\pa{1,0}}$ is a
$3$-independent generating set for $\FinCo{\omega_1}\times\FR{\omega_1}$.
That $X$ generates $\FinCo{\omega_1}\times\FR{\omega_1}$ is clear.  We use
proposition \ref{prop:a} to show that $X$ is $3$-independent.  Take any
$R\in\left[X\right]^{<\omega}$ and $\epsilon\in \leftexp R2$ such that $\prod_{x\in R}x^{\epsilon_x}=\pa{0,0}$.  Since there is no elementary product of elements of $\s{x_\alpha:\alpha<\omega_1}$ that is $0$, $\pa{1,0}\in R$ and $\epsilon_{\pa{1,0}}=1$.  Then there is a pair $a,b$ of elements in $R$ such that $\pi_2\pa a \perp \pi_2 b$ and $\epsilon_a=\epsilon_b=1$, so that $\s{\pa{1,0},a,b}\sub R$ and $\pa{1,0}\cdot a\cdot b =0$.  
\end{proof}

\section{Cardinal Function Results}\label{sec:cf}
Cellularity and independence have been considered earlier.  Here we give a few
results relating other cardinal functions to properties of $\perp$-graphs and
intersection graphs.  We will always assume that the graphs and algebras are infinite in
this section.



\begin{lemma}\label{finco_is_image}Let $A$ be $\omega$-free and $\omega
\leq\kappa=\left|A\right|$.  Then $B\defeq\FinCo\kappa$ is a homomorphic image of $A$.
\end{lemma}
\begin{proof}
Let $G$ be a set over which $A$ is $\omega$-free. Any bijective function $f:G\to \At B$ is $\omega$-preserving as all elements of $\At B$ are disjoint.  Since $A$ is $\omega$-free, $f$ extends to a homomorphism $\tilde f$ from $A$ to $B$.  Since the image of $f$ includes a set of generators, $\tilde f$ is surjective as well; that is, $B$ is a homomorphic image of $A$.
\end{proof}

The first use of this is that no infinite $\omega$-free \ba{} has the countable separation property. The \CSP{} is inherited by homomorphic images (5.27(c) in Koppelberg \cite{Handbook}), so if any infinite $\omega$-free \ba{} of size $\kappa$ 
has the \CSP{}, then by \ref{finco_is_image}, $\FinCo\kappa$ has the \CSP{}, which is a contradiction.  In particular, $\P\pa\omega /\mbox{fin}$ is not $\omega$-free.

We show that the spread of an $\omega$-free \ba{}  is  equal to its cardinality.

Theorem 13.1 of Monk \cite{CINV2} gives several equivalent definitions of spread, all of which have the same attainment properties; the relevant one to our purposes is the following.
$$s\pa A=\sup\left\{c\pa B : B\mbox{ is a homomorphic image of }A\right\}.$$

\begin{theorem}For $A$ $\omega$-free, $s\pa A=\left|A\right|$.  Furthermore, it is attained.\end{theorem}
\begin{proof}
From lemma \ref{finco_is_image}, $B=\FinCo{\left|A\right|}$ is a homomorphic image of $A$.  Since $c\pa B=\left|B\right|=\left|A\right|$, an element of the set in the above definition of $s\pa A$ is $\left|A\right|$.  Thus $s\pa A=\left|A\right|$ is attained.
\end{proof}


As they are greater than or equal to $s$, Inc, Irr, h-cof, hL, and hd are also equal to cardinality for $\omega$-free \bas{}.
Incomparability and irredundance are also attained by the $\omega$-free generating set.
This result also determines that $\left|\mbox{Id} A\right|=2^{\left|A\right|}$
as $2^{sA}\leq \left|\mbox{Id} A\right|$.  Then since $s$ is attained, $\left|\mbox{Sub}
A\right|=2^{\left|A\right|}$ as well.

%
%


The character of an $\omega$-free \ba{}  is also equal to cardinality.  Namely, at the bottom of page 183 in Monk \cite{CINV2}, it is shown that if $A$ is a homomorphic image of $B$, then $\chi\pa A \leq \chi \pa B$.  For $B$ $\omega$-free, let $A=\FinCo{\left|B\right|}$, so that $A$ is a homomorphic image of $B$ by lemma \ref{finco_is_image}, so we have that $\left|B\right|=\chi\pa A \leq \chi\pa B \leq \left|B\right|$.

\begin{theorem}If $A$ is infinite and  $\omega$-free, then  $\pi\pa{A}=\left|A\right|$.
\end{theorem}
Here $\pi$ is the density of $A$, the minimum of the cardinalities of dense subsets of $A$.
\begin{proof}
Take $H$ to be a set over which $A$ is $\omega$-free and let $D\sub A^+$ be dense.

For each $d\in D$, we can find a non-zero elementary product of elements of $H$ below $d$; write it as $\prod F_d \cdot \prod -G_d$ for finite disjoint $F_d,G_d\sub H$.

Now we show that $H=\bigcup_{d\in D} F_d$. Obviously $\bigcup_{d\in D} F_d\sub H$, so we need only show $H\sub\bigcup_{d\in D} F_d$. Choose an $h\in H$.  Since $D$ is dense, there is a $d\in D$ with $d\leq h$.  So $\prod F_d\cdot \prod -G_d \leq d\leq h$. Thus $\prod F_d \leq h+\sum G_d$, and since $H$ is $\omega$-independent, $h\in F_d$.

Since all the $F_d$ are finite, $\left|D\right|=\left|H\right|=\left|A\right|$
\end{proof}

We claim that the length (and therefore depth) of an $\omega$-free \ba{} is $\aleph_0$.  This
uses several preceding results.

\begin{theorem} If $A$ is $\omega$-free, then $A$ has no uncountable chain.
\end{theorem}
\begin{proof}Let $A$ be $\omega$-free over $G$.

Recall from theorem \ref{semigroup} that $A$ is a semigroup algebra over the set $H$ of finite products of elements of $G\cup\left\{0,1\right\}$.
 For $h\in H\sm\left\{0,1\right\}$, choose $g_1,\ldots,g_n\in G$ such that
$h=g_1\cdot\ldots\cdot g_n$ and set $h_G\defeq\left\{g_1,\ldots,g_n\right\}.$

Due to the result of Heindorf \cite{MR1155389}, if there is an uncountable chain in $A$,
there is an uncountable chain in $H$.  So by way of contradiction, we assume that there is an
uncountable chain $C\sub H$.  Without loss of generality, we may assume that
$0,1\not\in C$ so that every element of $C$ is a finite product of elements of
$G$. 

Let $C_G\defeq\left\{h_G:h\in C\right\}$.  We note that $$\bigcup C_G=\bigcup_{h\in
C}h_G\sub G$$ is the set of all elements of $G$ that are needed to generate the
elements of $C$, that is,

$C\sub \left<\bigcup C_G\right>$. so $C$ is a chain in that subalgebra of $A$
as well.

In order to reach a contradiction, we first show
that there are no finite subsets of $\bigcup
C_G$ with zero product.
Take $F\in\left[\bigcup C_G\right]^{<\omega}$.
Then for each $v\in F$, there is a $c_v\in C_G$ such that $v\in c_v$.  Note
that $\prod c_v\in C$ and $\prod c_v\leq v$.
Thus $\s{\prod c_v : v\in F}\sub C$, so $0\neq\prod\s{\prod c_v : v\in
F}\leq\prod F$, and hence $\prod F\neq 0$.

Thus $\bigcup C_G$ has no finite subset with zero product.  As $\bigcup C_G\sub G$ is
$\omega$-inde\-pendent, by lemma \ref{pindtoind}, it is independent.  Thus $\left<\bigcup
C_G\right>$ is free and hence has no uncountable chain, contradicting our original
assumption.
\end{proof}




\begin{theorem}\label{non_atom_endo}Let $A$ be $\omega$-free over $H$. Then $\left|\mathrm{End}\; A\right|= 2^{\left|A\right|}$.
\end{theorem}
\begin{proof}
For each $x\in H$, choose $y_x\in A$ such that $y_x<x$. For each $J\subset H$, define $f_J:H\rightarrow A$ as 
$$f_J\left(x\right)=\left\{\begin{array}{ll} y_x & x\in J \\ x & \mbox{otherwise.}\end{array}\right.$$
$f_J$ is $2$-preserving and extends to an endomorphism.  So we have exhibited $2^{\left|A\right|}$ endomorphisms.
\end{proof}

%% file: OtherBAs.tex
\section{Maximal $n$-independence number}
We can look at $n$-independent sets in \bas{} that aren't $n$-free.  The natural thing to do is introduce a cardinal function, $\PInd$, that measures the supremum of the cardinalities of those sets. Since $\PInd$ is a regular sup-function, we can define a spectrum function and a maximal $n$-independence number of a \ba{} in the standard way.
\begin{definition}Let $1\leq n\leq \omega$.
$$\J_{nsp}\pa A \defeq\left\{\left|X\right|: X\;\mbox{is a
maximal}\;n\!\mbox{-independent subset of}\; A\right\}$$
$$\J_n\pa A \defeq \min\pa{\J_{nsp}\pa A}$$
\end{definition}
This could be written as $\PInd_{mm}$ according to the notation of
Monk \cite{CINV2}.  Note that $\J_1=\mathfrak i$ where $\mathfrak i$ is the
minimal independence number as seen in Monk \cite{ContCard}.

This is defined for every \ba{}; from the definition it is easily seen that the
union of a chain of $n$-independent sets is $n$-independent, so Zorn's
lemma shows that there are maximal $n$-independent sets.  $\J_n\pa{A}$ is infinite for all $n\leq\omega$ if $A$ is atomless (shown in lemma \ref{mpi_sum_1}), and has value $1$ if $A$ has an atom.

\begin{proposition} For all $n$ with $1\leq n\leq\omega$, if $A$ has an atom, then $\J_n\pa A=1$.\end{proposition}
\begin{proof}
If $a$ is an atom of $A$, then we claim that $\s{-a}$ is a maximal $n$-independent subset of $A^+$.  That $\s{-a}$ is $n$-independent is clear as any singleton other than $\s 0$ and $\s1$ is independent. 

Let $x\in A^+\sm\s{-a}$, we show that $\s{-a,x}$ is not $n$-independent.
There are two cases.

If $a\leq x$, then $1=a+-a\leq x + -a$, so that $\pa{\perp1}$ fails.

If $a\leq -x$, then $x\leq -a$, so that $0\neq\prod\s x\leq\sum\s{-a}$, but $\s
x\cap\s{-a}=\emptyset$, so that $\pa{\perp3}$ fails.
\end{proof}

\begin{lemma}\label{mpi_sum_1} Let $B$ be a \ba{}, $2\leq n\leq \omega$, and $H\sub B^+$ be
$n$-independent.  If $H$ is maximal among $n$-independent subsets
of $B^+$, then $H$ is infinite and $\sum H =1$ or $H$ is finite and $-\sum H$ is an atom.\end{lemma}
\begin{proof}We prove the contrapositive.
First, the case that $H$ is infinite.
Let $H\sub B^+$ be $n$-independent and have $b<1$ as an upper bound. 
We show that $H\cup\s{-b} $ is $n$-independent:

Note that $-b\notin H$, as $-b\not\leq b$.  Now we will apply proposition
\ref{prop:a}.  So, assume that $R\in\subs{H\cup\s{-b}}{<\omega},
\epsilon\in\leftexp R2$, and $\prod_{x\in R}x^{\epsilon_x}=0$.  If
$-b\notin R$, the conclusion follows since $H$ is $n$-independent.  So suppose
that $-b\in R$.  Let $R'\defeq R\sm\s{-b}$.  Then we have two cases:
\begin{enumerate}
\item[Case 1.] $\epsilon_{-b}=1$.  If there is an $x\in R'$ such that
$\epsilon_x=1$, then $x\leq b$ and so $x\cdot -b =0$ as desired.  So assume
that $\epsilon\left[R'\right]=\s0$.  Then $-b\leq\sum_{x\in R'}x\leq b$, which
is a contradiction.
\item[Case 2.] $\epsilon_{-b}=0$.  If $\epsilon_x=1$ for some $x\in R'$, then 
$$0=\prod_{y\in R} y^{\epsilon_y}=\prod_{y\in R'}y^{\epsilon_y}\cdot
b=\prod_{y\in R'}y^{\epsilon_y}$$
and the $n$-independence of $H$ gives the result.  So assume that
$\epsilon\br{R'}=\s0$.  Then $b\leq\sum R'\leq b$, so $b=\sum R'$.  Then
$b\cdot \prod_{x\in R'}-x=0$, contradicting the $n$-independence of $H$.
\end{enumerate}
So we have that if $H$ is infinite and maximal $n$-independent, it has no upper bound
other than $1$, so $\sum H =1$.

Now we consider the case that $H$ is finite.
 If $-\sum H$ is not an atom, let $0<a<-\sum H$, then we claim that $H\cup\s a$
is $n$-independent. Again we use proposition \ref{prop:a}.  Assume that
$R\in\subs{H\cup\s a}{<\omega}, \epsilon\in\leftexp R2$, and $\prod_{x\in
R}x^{\epsilon_x}=0$. \Wolog, $a\in R$. 
\begin{enumerate}
\item[Case 1.]$\epsilon_a=1$.  If $\epsilon_x=1$ for some $x\in R\sm\s a$, then
$a\cdot x\leq a\cdot\sum H=0$, as desired.  Otherwise $$a\leq \sum\pa{R\sm\s
a}\leq \sum H$$ and so $a=0$, contradiction.
\item[Case 2.]$\epsilon_a=0$.  If $\epsilon_x=1$ for some $x\in R\sm\s a$, then
$a\cdot x=0$, hence $x\leq-a$, and then $$\prod_{y\in
R}y^{\epsilon_y}=\prod\s{y^{\epsilon_y}:y\in R\sm\s a}$$ and the conclusion
follows.  Otherwise $$-a\leq\sum\pa{R\sm\s a}\leq\sum H,$$ so $-\sum H\leq a$,
contradicting $a<-\sum H$.
\end{enumerate}

\end{proof}

The converse of lemma \ref{mpi_sum_1} does not hold.  An example due to Monk is in  $\FR{X\cup Y}$ where $X\cap Y=\emptyset$ and 
 $\left|X\right|=\left|Y\right|=\kappa\geq\omega$.
  $X$ is independent, is not maximal for $2$-independence, and has sum 1.  Here $\sum X=1$ is the only non-trivial part--by way of contradiction, let $b$ be a non-1 upper bound for $X$.  Then $-b$ has the property that $x\cdot -b =0$ for all $x\in X$, so let $a$ be a elementary product of elements of $X\cup Y$ where $a\leq -b$.  Take some $x\in X$ that does not occur in that elementary product.  Then since $X\cup Y$ is independent, $a\cdot x \neq 0$, but since $a\leq -b$, $a\cdot x =0$.

\begin{theorem}For $B$ atomless, and $2\leq n\leq \omega$, $\mathfrak p\pa B \leq \J_n \pa B$.\end{theorem}
Here $\mathfrak p\pa B$ is the pseudo-intersection number, defined in
Monk \cite{ContCard} as 
$$\mathfrak p \pa A \defeq \min\s{\left|Y\right|:Y\sub A\;\mbox{and}\;\sum Y = 1 \;\mbox{and}\; \sum Y'\neq 1\;
\mbox{for every finite}\; Y'\sub Y}.$$
\begin{proof}
Since $B$ is atomless, a maximal $n$-independent set $Y$ has $\sum Y=1$,
and by \pc1, if $Y'\sub Y$ is finite, $\sum Y'\neq 1$.  That is, the
maximal $n$-independent sets are included among the $Y$ in the definition
of $\mathfrak p$. 
\end{proof}

We do not know if strict inequality is possible.

\begin{corollary} For all $n$ with $1\leq n\leq \omega$, $\J_n\pa{\PoF}\geq\aleph_1$\end{corollary}
\begin{proof}$\aleph_1\leq \mathfrak p \pa{\PoF} \leq \J_n\pa{\PoF}$
\end{proof}

We also recall that under Martin's Axiom, $\mathfrak p \pa{\PoF} =
\beth_1$, so the same is true of $\J_n$.

\begin{proposition} Any $B$ with 
the strong \CSP{} has, for all $2\leq n\leq \omega$,  $\J_n\pa B\geq \aleph_1$.\end{proposition}
\begin{proof}Such a $B$ is atomless, so let $H\sub B^+$ be
$n$-independent and countably infinite, that is $H=\left<h_i:i\in \omega\right>$. Then
let $c_m\defeq\sum_{i\leq m} h_i$. Each $c_m$ is a finite sum of elements of $H$,
thus by $\pa{\perp1}$, $c_m<1$. Then $C\defeq\s{c_i:i\in \omega}$ is a countable
chain in $B\sm \s1$, so by the strong \CSP{}, there is a $b\in B$ such that
$c_i\leq b <1$ for all $i\in \omega$.  Then as $h_i\leq c_i$, $h_i\leq b$ for
all $i\in\omega$ as well, that is, $b$ is an upper bound for $H$.  Thus by
lemma \ref{mpi_sum_1}, $H$ is not maximal.
\end{proof}

In addition, we show that maximal $n$-independent sets lead to weakly dense sets. 

We use the notation $-X=\s{-x : x\in X}$ frequently in the sequel.

\begin{theorem}Let $1\leq n\leq\omega$.  If $X\sub A$ is maximal $n$-independent in $A$, then the set $Y$ of nonzero elementary products of elements of $X$ is weakly dense in $A$.
\end{theorem}
Recall that $Y$ is weakly dense in $A$ if and only if $Y\sub A^+$ and for every $a\in A^+$, there is a $y\in Y$ such that $y\leq a$ or $y\leq -a$.
\begin{proof}
If $a\in X$, this is trivial, so we may assume that $a\notin X$ and hence $X\cup\s a$ is not $n$-independent.

By proposition \ref{prop:a}, there exist $R\in\subs{X\cup\s a}{<\omega}$ and
$\epsilon\in \leftexp R2$ such that $\prod_{x\in R}x^{\epsilon_x}=0$ while for
every $R'\in\subs R{\leq n}$, if $\epsilon\br{R'}\sub \s 1$ then $\prod
R'\neq0$.  This last implication holds for every $R'\in\subs{R\sm\s a}{\leq
n}$, and so $\prod\s{x^{\epsilon_x}:x\in R\sm\s a}\neq 0$ since $X$ is
$n$-independent.  But $\prod\s{x^{\epsilon_x}:x\in R\sm\s a}\leq a$ or $\leq
-a$, as desired.
\end{proof}

\begin{corollary}
If $A$ is atomless and $1\leq n\leq \omega$, then $\mathfrak r\pa A\leq \J_n\pa A$.
\end{corollary}
Recall the definition of the reaping number: $$\mathfrak r\pa A\defeq\min\s{\left|X\right| : X\mbox{ is weakly dense in } A}$$
\begin{proof}
Since $A$ is atomless, all maximal $n$-independent sets are infinite, and thus there is a set of size $\J_n\pa A$ weakly dense in $A$.
\end{proof}
We do not know if strict inequality is possible.

We do not currently have any results for the behavior of $\J_n$ on any type of
product or its relationship to $\mathfrak u$.

We show the consistency of $\J_n\pa{\P\pa \omega /fin} < \beth_1$ for $1\leq n\leq\omega$.  The argument is similar to
exercises (A12) and (A13) in chapter VIII of Kunen \cite{KunenSet}; the main lemma
follows.
 
\begin{lemma}\label{con_j_lemma}
Let $M$ be a countable transitive model of $ZFC$ and $1\leq k\leq\omega$.  For a subset $a$ of
$\omega$, let $\left[a\right]$ denote its equivalence class in $\P\pa\omega /
fin$.  Suppose that $\kappa$ is an
infinite cardinal and $\left<a_i : i<\kappa\right>$ is a system of infinite
subsets of $\omega$ such that $\left<\left[a_i\right]:i<\kappa\right>$ is
$k$-independent in $\P\pa\omega / fin$.  Then there is a generic extension
$M\left[G\right]$ of $M$ using a ccc partial order such that in
$M\left[G\right]$ there is a $d\sub \omega$ with the following properties:
\begin{enumerate}
\item $\left<\left[a_i\right]:i<\kappa\right>^\frown \left<\left[\omega\sm
d\right]\right>$ is $k$-independent.

\item If $$x\in \pa{\P\pa\omega \cap M}\sm\pa{\s{a_i : i<\kappa}\cup\s{\omega\sm
d}},$$ then $$\left<\left[a_i\right]:i<\kappa \right>^\frown\left<\left[\omega\sm
d\right],\left[x\right]\right>$$ is not $k$-independent.
\end{enumerate} 
\end{lemma}
\begin{proof}We work within $M$ here.

Let $B$ be the $k$-independent subalgebra of \PoF{} generated by $\s{\left[a_i\right]:i<\kappa}$.  By Sikorski's extension criterion, let $f$ be a homomorphism from \\$\left<\s{a_i : i<\kappa}\cup \s{\s m : m\in \omega}\right>$ to $\overline B$ such that $f\pa{a_i}=\left[a_i\right]$ and $f\pa{\s m}=0$.  
Then let $h:\P\pa\omega \fct \overline B$ be a homomorphic extension of $f$ as given by Sikorski's extension theorem.

Let $ P\defeq\s{\pa{b,y}:b\in\ker\pa h \mbox{ and }
y\in\left[\omega\right]^{<\omega}}$ with the partial order given by
$\pa{b,y}\leq\pa{b',y'}$ if and only if $b\supseteq b'$, $y\supseteq y'$ and $y\cap b'\sub
y'$.  This is a ccc partial order.  Let $G$ be a $ P$-generic filter
over $M$, and let $d\defeq\bigcup_{\pa{b,y}\in G}y$.

We now have several claims that combine to prove the lemma.

\begin{itemize}

\item[Claim 1.] If $R$ is a finite subset of $\kappa$ and $\epsilon\in \leftexp R 2$ is such that $\bigcap_{\stackrel{i\in R}{\epsilon_i=1}}a_i$ is infinite, then $\bigcap_{i\in R}a_i^{\epsilon_i}\cap d$ is infinite.

Let $R$ and $\epsilon$ be as given, then for each $n\in \omega$, let 
$$E_n\defeq \s{\pa{b,y}\in P : \exists m>n\left[m\in\bigcap_{i\in R}a_i^{\epsilon_i}\cap y\right]}.$$

First, we show that each $E_n$ is dense.  Take $\pa{b,y}\in  P$.  Then
$c\defeq\pa{\bigcap_{i\in R}\pa{a_i^{\epsilon_i}}}\sm b$ is infinite; if
not, then $c$ is finite (thus in $\ker\pa h$, as is $b$)and $\bigcap_{i\in R}a_i^{\epsilon_i}\sub b\cup c$.  Applying $h$ to both sides gives $\prod_{i\in R}\left[a_i\right]^{\epsilon_i}=0$, which is a contradiction of proposition \ref{prop:a}.  So we choose an $m\in c\sm y$ such that $m>n$; then $\pa{b,y\cup\s m}\leq\pa{b,y}$ and  $\pa{b,y\cup\s m}\in E_n$, showing that $E_n$ is dense. 
This shows the claim, as for each $n\in\omega$, $E_n\cap G\neq \emptyset$, so that we have an integer larger than $n$ in $\bigcap_{i\in R}a_i\cap d$.

\item[Claim 2.]  If $R$ is a finite subset of $\kappa$ and $\epsilon\in \leftexp R 2$ such that $\bigcap_{\stackrel{i\in R}{\epsilon_i=1}}a_i$ is infinite, then $\bigcap_{i\in R}a_i^{\epsilon_i}\sm d$ is infinite. 

Let $R$ and $\epsilon$ be as given, then for each $n\in \omega$, let 
$$D_n \defeq\s{\pa{b,y}\in P:\exists m>n\left[m\in\bigcap_{i\in R}a_i^{\epsilon_i}\cap b\sm y\right]}. $$

To show that $D_n$ is dense, take any $\pa{b,y}\in P$.  Since $\bigcap_{i\in
R}a_i^{\epsilon_i}$ is infinite from proposition \ref{prop:a}, it follows that
we may choose $m>n$ such that $m\in \bigcap_{i\in R}a_i^{\epsilon_i}\sm y$.
Then $\pa{b\cup\s m , y}\leq\pa{b,y}$ and $\pa{b\cup\s m,y}\in D_n$, as
desired.

Take some $\pa{b,y}\in D_n\cap G$.  Then there is an $m>n$ such that $m\notin d$ (thus proving the claim). In fact, choose $m>n$ such that $m\in\bigcap_{i\in R}a_i^{\epsilon_i}\cap b \sm y$. We claim that $m\not\in d$. Suppose that $m\in d$;
 then we have a $\pa{c,z}\in G$ with $m\in z$ and $\pa{e,w}\in G$ that is a common extension of $\pa{b,y}$ and $\pa{c,z}$.  Then $m\in w\cap b\sm y$, contradicting that $\pa{e,w}\leq\pa{b,y}$.

\item[Claim 3.]
$\left<\left[a_i\right]:i<\kappa\right>^\frown\left<\left[\omega\sm
d\right]\right>$ is $k$-independent.

Suppose that $R\in\subs\kappa{<\omega},\epsilon\in\leftexp R2, \delta\in 2$,
and $\prod_{i\in R}\br{a_i}^{\epsilon_i}\cdot\br{\omega\sm d}^\delta =0$.  By
claims 1 and 2 (depending on $\delta$), $\prod_{\stackrel{i\in
R}{\epsilon_i=1}}\br{a_i}^{\epsilon_i}=0$.  Since
$\left<\br{a_i}:i<\kappa\right>$ is $k$-independent, there is a subset
$R'\sub\s{i\in R:\epsilon_i=1}$ of size at most $k$ such that $\prod_{i\in
R'}\br{a_i}=0$, as desired.

\item[Claim 4.] If $b\in \ker\pa h$, then $b\cap d$ is finite.

$\s{\pa{c,y}\in P:b\sub c}$ is dense in $P$, so that there is a $\pa{c,y}\in G$
such that $b\sub c$.  We show $b\cap d \sub y$ and thus is finite.  Let $m\in
b\cap d$ and choose an $\pa{e,z}\in G$ such that $m\in z$.  Let $\pa{r,w}\in G$
be a common extension of $\pa{e,z}$ and $\pa{c,y}$; then (recalling the
definition of the order) $m\in w\cap c\sub y$.

\item[Claim 5.]If $$x\in
\pa{\P\pa\omega\cap M}\sm\pa{\s{a_i:i<\kappa}\cup\s{\omega\sm d}},$$ then
$$s\defeq\left<\left[a_i\right]:i<\kappa\right>^\frown\left<\left[\omega\sm
d\right],\left[x\right]\right>$$ is not $k$-independent.

We have two cases here.  The slightly easier is if $x\in\ker\pa h$; then by claim 4,
$x\cap d$ is finite, so that $\left[x\right]\leq\left[\omega\sm d\right]$,
 causing $s$ to fail to even be ideal-independent.  If
$x\notin\ker\pa h$, then there is a $b\in B$ with
$0<b\leq h\pa x$.  Since $B$ is $k$-freely generated by
$\left<\left[a_i\right]:i<\kappa\right>$, we may take $b$ to be a elementary product
of elements of $\left<\left[a_i\right]:i<\kappa\right>$.  Then
$b=\left[c\right]$, where $c=\bigcap_{i\in R}a_i^{\epsilon_i}$ is infinite.  Then $c\sm x\in \ker\pa h$.  By claim 4, this gives $\prod_{i\in R}\left[a_i\right]^{\epsilon_i}\cdot -\left[x\right]\cdot\left[d\right]=0$, contradicting proposition \ref{prop:a} for $s$.

\end{itemize}
\end{proof}

\begin{theorem}For each $1\leq k\leq\omega$, it is consistent with $\beth_1>\aleph_1$ that \\$\J_k\pa{\PoF}=\aleph_1$.\end{theorem}
\begin{proof}We begin with a countable transitive model $M$ of $ZFC + \beth_1>\aleph_1$, then iterate the construction of lemma \ref{con_j_lemma} $\omega_1$ times as in lemma 5.14 of chapter VIII of Kunen \cite{KunenSet}. This results in a model of $ZFC + \beth_1>\aleph_1 + \J_k\pa{\PoF}=\aleph_1$.
\end{proof}

This shows that $\J_k\pa{\PoF}=\beth_1$ is independent of $ZFC$.